\documentclass[12pt]{article}
 \usepackage{amsmath,amsfonts,theorem}
\input epsf
 \newcommand{\ntag}[1]{} 
\numberwithin{equation}{section}
 
 \newtheorem{prop}{Proposition}

 \newtheorem{cor}{Corollary}

 \theorembodyfont{\rmfamily}
 \newtheorem{dfn}{Definition}

 \newcommand{\qed}{\ifhmode\unskip\nobreak\fi\quad\ensuremath\square}

 \newcommand{\sA}{\mathcal A} 
 \newcommand{\sB}{\mathcal B} 
 \newcommand{\sG}{\mathcal G} 
 
 \newcommand{\sH}{\mathcal H}

 \newcommand{\sM}{\mathcal M}

 \newcommand{\al}{\alpha}

 \newcommand{\ga}{\gamma}
 
 \newcommand{\om}{\omega}
 
 \newcommand{\De}{\Delta}
 \newcommand{\Ga}{\Gamma}
 
 \newcommand{\Om}{\Omega}

 \newcommand{\PP}{\mathbb P}
 \newcommand{\C}{\mathbb C}

 \newcommand{\R}{\mathbb R}
 \newcommand{\Z}{\mathbb Z}

 \renewcommand{\Im}{\operatorname{Im}}

\begin{document}

 \title{Maslov class of  lagrangian
 embedding to Kahler manifold}

  \author{N.A. Tyurin}

\date{}

\maketitle

\begin{center}
 {\em BLThP JINR (Dubna) and MSU of RT}
 \end{center}
 \bigskip

\section*{Introduction}

The development of symplectic geometry during the last twenty years
of the previous century and the first years of the present one is
originated mostly on the tendency to transport or reformulate numerous results,
derived in classical mechanics for symplectic vector spaces, to
the case of compact symplectic manifolds. Although various
problems of special interest are formulated
 for  the compact case only (f.e., the Arnold conjecture),
it is natural to move from simple cases   to complicated ones, from the flat case to
the compact one. It is  specifically  important for such  problems as the
problem of quantization of classical mechanical systems, since the physical
background of the problem requires the compatibility of the results
for the flat and the compact cases. Namely (see, f.e.,  \cite{Wood}),
a quantization scheme can be acknowledged as a working one if
its reduction to the case of the symplectic vector space agrees
with the canonical quantization, a dogma of modern
theoretical physics.

One of the essential ingredients of the symplectic geometry of
the vector space is the Maslov class, induced by  lagrangian embedding.
Recall (see \cite{ArGiv}), if $(V, \om)$ --- a symplectic vector
space with a constant  symplectic form $\om$ and $S \subset V$
is a lagrangian submanifold, smoothly embedded\footnote{in the present
paper we consider  the case of  smooth embeddings only, although
our constructions  with some corrections can be generalized
to the case of lagrangian immersions} to $V$, then this embedding
 induces a class $m_S \in H^1(S, \Z)$,
 which is integer valued and it implies its most important
 property --- the invariance under continuous deformations preserving
 the lagrangian condition. It follows that the class $m_S$
 is stable with respect to the hamiltonian deformations of $(V, \om)$
 and is equivariant under the symplectomorphisms
 of our vector space. Most important applications
 of the Maslov class are the application for the minimality
 problem, solved in \cite{Mor}, and the realization of the class
 as a correction term in the quantization method of Karasev - Maslov
 which is often called   the quasi classical quantization, see \cite{KarMas}.
 Specifically the second application attracted our attention to
 the possibility of generalizations of the Maslov class to
 the case of  lagrangian
 embeddings satisfied some special conditions with respect
 to the Levi - Civita connection to a Kahler manifold,
 to exploit it in constructions of algebraic lagrangian
 geometry (see \cite{GT}) and ALG(a) - quantization (see \cite{jT1}).

The starting point of the present work
was the attentive reading of the paper \cite{T1}, were the Maslov
class is related to the abelian connection on the restricted
to a lagrangian submanifold anticanoncial line bundle induced by
the Levi - Civita connection. In other words, the present paper
is just the explanation of three lines of the paper \cite{T1}.
Usually one deduces a generalized Maslov class from the mean curvature
of a lagrangian submanifold (see \cite{Oh1}, \cite{Ars})
which is a real 1 -form, and in this circumstance the minimality
problem is solved automatically since a lagrangian submanifold has
the minimal riemannian volume in the class of local deformations if and only if
its mean curvature vanishes (see \cite{Br}, \cite{Oh1}).
However even if the mean curvature of a lagrangian submanifold
is a closed form it is rather difficult to show that
its cohomology class is integer valued (but it is done
in several cases, see \cite{Oh1}).  In the papers  \cite{T1}
and \cite{T2}  in the special situation of the Kahler - Einstein
manifolds  some universal Maslov class is defined via connections,
and this automatically implies that it is integer valued, but
now it is necessary to relate this class to the mean curvature.
The main technical result of the present paper is exactly
to establish the desired relationship in the general case.
Namely, Proposition 9 states and proves that for any
oriented lagrangian submanifold $S \subset M$
of a Kahler manifold $M$ the mean curvature
$\al_H$ up to scaling coincides with the 1 - form, which represents
the determinant Levi - Civita connection on the anticanonical line
bundle with respect to the canonical flat structure.
Since the curvature of the determinant Levi - Civita connection
multiplied by $\frac{1}{2 \pi i}$ coincides with the Ricci
form $\rho$ of the Kahler manifold,  from this fact one gets the following
 identity
$$
d \al_H = \rho|_S,
$$
which is already known from  \cite{Daz}.

The construction presented in section 1
endows a lagrangian submanifold of a Kahler manifold, such that
the restriction to it of the determinant Levi - Civita connection
is trivial and admits covariant constant sections, with a map
to $S^1$, which is called  {\it the phase};
in the case when the Maslov class is trivial
the phase is represented by a smooth function
$$
\phi_S: S \to U(1),
$$
which is defined up to a constant phase scaling (the logarithmic
derivative of the "total" phase gives a closed 1- form, whose cohomology
class coincides with the Maslov class $m_S$). For the case
of the Kahler - Einstein manifolds when the Ricci tensor
is proportional to the symplectic form,
$$
\rho = k \om,
$$
it makes possible to define
some half weighting rules for Bohr - Sommerfeld cycles,
which are studied in  \cite{GT}. It is the main
result of the paper, and we hope that it will be exploited
in algebraic lagrangian geometry, constructed by A. Tyurin
and A. Gorodentsev in \cite{GT}.

At the time he was a student of V.I. Arnold
the author could not present the adequate diploma work and
now he hopes that this paper can be regarded as such a work
at the same time expressing author's gratitude to the former
scientific advisor. The present work couldn't be done without
help and remarks of A. Gorodentsev. The author would
like to express his deep gratitude to D. Orlov, D. Auroux
and S. Nemirovskiy. During the work on this paper the author
was supported by the Max Planck Institute for Mathematics
(Bonn, Germany) and European Center for Nuclear Research (Geneva,
Switzerland) and I would like to cordial thank
all staff of both the institutions. The work was partially supported
by RFBR (grants NN 05 - 01 - 01086, 05 - 01- 00455).

\section{The phase and the Maslov class}

The definition of the Maslov class in the classical case
 --- for a lagrangian embedding to a symplectic vector space ---
 can be generalized to the case of any Kahler manifold for lagrangian
 submanifolds which possess special properties with respect to
 the determinant Levi - Civita connection. At the same time
 in the most generic case the Maslov class is not  invariant
 under hamiltonian deformations, however restricting
 the construction to the class of Kahler - Einstein manifolds
 the Maslov class is an invariant of isodrastic deformations, and
 for even more specified case --- the case of Calabi - Yau manifolds ---
 it is invariant  for all lagrangian deformations.

The basic construction is as follows. Let $(M, \om, I)$ be any Kahler
manifold which is understood as a symplectic manifold $(M, \om)$
equipped with an integrable complex structure $I$ compatible with $\om$.
Equivalently the manifold can be understood as a symplectic manifold equipped
with a compatible riemannian metric $G$ with holonomy $U(n)$, where $\dim_{\R} M = 2n$.
These representations are equivalent. In absolutely general non integrable
case a riemannian metric induces the Levi - Civita connection
on the tangent bundle $TM$, and in the integrable case this connection
commutes with the operator $I$. This means that the Levi - Civita
connection reduces to a hermitian connection on the holomorphic
tangent bundle  $T^{1,0} M$. Its determinant  $\det T^{1,0} M$
is a complex line bundle which is called the anticanonical bundle
and is denoted as
$K^{-1}_M$; it carries a canonical hermitian structure
and an abelian connection  $A_{LC}$, induced by the total Levi - Civita
connection. To avoid confusions in the rest of the paper we will call
 $A_{LC}$ on $K^{-1}_M$ the determinant Levi - Civita connection
 (since it really is the determinant of the total Levi - Civita connection).
The curvature of this connection,
$F(A_{LC})$, up to scaling coincides with the Ricci
form $\rho$ of the Kahler structure on $M$ (all details can be found
in \cite{GrHar}).

Further, let $S \subset M$ be a smooth oriented lagrangian
submanifold: it means that  $\dim S = \frac{1}{2} \dim M = n$
and
$$
\om|_S \equiv 0.
$$
Then it's not hard to see that the restriction of the anticanonical bundle
to $S$ has the form
$$
K^{-1}|_S = \det TS \otimes \C.
\eqno (1)
$$
In particular, it implies that the restriction of the (anti) canonical
class to any oriented lagrangian submanifold is topologically
trivial. This is true as well in the case of non integrable complex
structures. Moreover, the compatibility condition for  $\om$ and
 $I$ implies that one has even more generic fact:

\begin{prop} The holomorphic tangent bundle $T^{1,0} M$ being restricted to
$S$ has the form
$$
T^{1,0} M |_S = TS \otimes \C;
$$
moreover, the hermitian structure on $TM|_S$, defined by the hermitian
triple $(\om, I, G)$, is the complexification of the special
orthogonal structure, defined by the riemannian metric $G$ on $TS$.
\end{prop}

Indeed, the proof\footnote{suggested to the author
by Dennis Auroux, and I want to express the cordial thanks
for this hint} just requires to apply the Darboux - Weinstein theorem,
 \cite{Wein}, which describes neighborhoods of lagrangian submanifolds.
 The restriction of the real tangent bundle $TM|_S$
 can be filled by the deformation family of $TS$ under the
 action of operators
 $$
 \cos \phi \cdot {\rm Id} +  \sin \phi \cdot {\rm
 I}, \quad 0 \leq \phi < \pi.
 $$

 We will use Proposition 1 in section 3, and now it is sufficient
 for our aims to consider its reduction to the  determinants.
 Namely, the identification (1) is completed by the considerations
 of the canonical hermitian structure  on $K^{-1}|_S$
 and the special orthogonal structure on  $\det TS$ which are compatible.
 This means that the space of hermitian connections  $\sA_h(K^{-1}|_S)$
 containes two distinguished elements: the determinant Levi -
 Civita connection and the trivial connection, defined by
 the complexification of the determinant part of the restricted to $S$
 Levi - Civita connection. Indeed, let's restrict the metric
 $G$ to $S$ and consider a trivialization of $\det TS$, defined by
 an appropriate polivector field of constant length (it is dual
 to a volume form $d \mu$). Simultaneously this trivialization
 in view of the isomorphism  (1) defines a trivialization of $K^{-1}|_S$,
 and the corresponding trivial connection is denoted as $A_0$.

Now, suppose that an oriented lagrangian embedding $S \subset M$ satisfies
a strong condition  ---  the determinant connection $A_{LC}$
is {\it flat } and has  {\it trivial periods} on $H_1(S, \Z)$. It's well known
(see, f.e., \cite{DK}) that the abelian connection space has very simple
structure; in particular in our case the connections $A_{LC}$ and $A_0$
are  {\it gauge equivalent}. They are related by a gauge transformation
 $$
 g_S = g(A_{LC}, A_0) \in {\rm Map} (S, U(1))/ U(1) = {\rm Map} (S, S^1).
\eqno (2)
 $$

From this one gets two definitions for an oriented lagrangian embedding
with trivial restriction of the determinant Levi - Civita connection:
\begin{dfn} The gauge transformation  $g_S$ (2) is called the phase
of the lagrangian submanifold $S \subset M$.
\end{dfn}

\begin{dfn} The  Maslov class of the lagrangian embedding
 $S \subset M$ to the Kahler manifold $(M, \om, I)$ is
 the integer valued cohomology class on $S$
$$
m_S = g^*_S h,
$$
(see (2)), where $h \in H^1(S^1, \Z)$ is the generator.
\end{dfn}

It is convenient to represent the phase (at least locally in the case of nontrivial
Maslov class) in the form

$$
e^{i \phi_S} : S \to U(1),
$$
instead of the gauge transformation (2), that is described by a real function
$\phi_S$, defined up to an additive constant. By misuse of language
we will sometimes call the function by the same name "phase" if it does
 not lead to confusion.

In general case, the logarithmic derivative of $e^{i \phi_S}$ gives a closed
1- form whose cohomology class coincides with the Maslov class $m_S \in H^1(S, \Z)$.

One finishes the construction of the Maslov class for a lagrangian
embedding with trivial restriction of the determinant Levi - Civita
connection with the following

\begin{prop} The definitions of the phase $g_S$
and the  Maslov class $m_S$ for a lagrangian embedding with trivial
restriction of the determinant Levi - Civita connection
are correct.
\end{prop}

Indeed, the construction doesn't
need any additional data and hence is universal.

In general for a specified problem the phase can be defined for any lagrangian
embedding however for this one needs to introduce some special rules
and thus one loses the universality. For example, one can define
an abelian connection $A^0_{LC}$, which is the result
of the orthogonal projection of the determinant Levi - Civita
connection $A_{LC}$ to the gauge orbit $\sG_h(A_0) \subset \sA_h(K^{-1}|_S)$
(all details of the theory of connections and curvatures
as well as the Hodge theory of harmonic forms can be found
in \cite{DK}).

Consider the curvature of the connection $A_{LC}$:
$$
F(A_{LC}) = F_{LC} \in i \Om^2_S.
$$
From the Chern - Weyl theory we know that the curvature form $F_{LC}$,
multiplied by $\frac{i}{2 \pi}$, is a closed real form which represents
the cohomology class $c_1(K^{-1}|_S) \in H^2(S, \Z)$. But this class
is trivial (see (1)), thus  $\frac{i}{2 \pi} F_{LC}$ is an exact form.
Choose the minimal (with respect to the riemannian metric) correction term
 $\De_1$ to connection $A_{LC}$ such that $A_{LC} - \frac{i}{2 \pi} \De_1$
 is flat. According to the Hodge theory such a form exists and is unique:
 from the set of forms which satisfy
$$
d \De_1 = \frac{i}{2\pi} F_{LC},
$$
we choose the correction term imposing additionally that
$\De_1 \bot \ker d$. Further, the connection
$$
A^1_{LC} = A_{LC} - \frac{i}{2\pi} \De_1
$$
is already flat, and one can consider its character on the fundamental
group $\pi_1(S)$. But the connections we study are abelian hence the character
of the connection $A^1_{LC}$ on the commutant of the group
 $\pi_1(S)$ is trivial, and it descends
to the periods of the connection on the lattice
$$
H_1(S, \Z) = \pi_1(S) / [\pi_1(S), \pi_1(S)].
$$
 The periods of the connection $A^1_{LC}$ on $H_1(S, \Z)$ are
 uniquely defined by the corresponding class from $H^1(S, \R)$,
 and there are such harmonic forms $\De_2 \in \sH^1_S$, that
 connection
$$
A^0_{LC} = A_{LC} - \frac{i}{2\pi} \De_1 - \frac{i}{2\pi} \De_2
$$
has trivial periods on the lattice $H_1(S, \Z)$ and hence the trivial
character on $\pi_1(S)$. In some cases one can choose  unique correction
form, distinguished by the minimality condition with respect to the riemannian
norm, and such a minimal form is unique only if the periods of the connection
$A^1_{LC}$ on the primitive elements of $H_1(S, \Z)$ don't have half integer
values. If this condition holds then we have
 uniquely defined correction term $\De_2$,
correctly defined trivial connection $A^0_{LC}$
and again we define the gauge transformation $g_S$, which
transport  $A^0_{LC}$ to $A_0$. The exception is given by
some "boundary" case when at least one period is half integer:
$$
{\rm Mon} A^1_{LC} (\ga) = -1, \quad \ga \in H_1(S, \Z),
$$
In this case it is necessary either to define half integer
Maslov classes or to introduce some additional data
to choose $\De_2$. If all the periods are not
half integer then the correction term $\De_2$ is uniquely defined
and the projection $A^0_{LC}$ is correctly defined.

Anyway, for an oriented lagrangian embedding $S \subset M$
to a Kahler manifold $(M, \om, I)$, we can decompose the connection
 $A_{LC}$ into three parts:

--- the part which responses to the curvature $F_{LC}$ (that is to
the restriction of the Ricci form $\rho$), denoted as $\De_1$,

--- the part which corresponds to the periods of $A_{LC}^1$, that
is a point on the Jacobian
$$
H^1(S, \R)/ H^1(S, \Z),
$$
 denoted as $\De_2$,

--- the phase, a gauge transformation, denoted as $g_S$.

Let's emphasize again that the choice
of $\De_2$ is not universal and can be done universally
just in certain specified cases. For example, if we consider
lagrangian embeddings with bounded periods that is the embeddings
for which the periods of the flat part $A^1_{LC}$
on the primitive elements of the lattice $H_1(S, \Z)$
strictly less than $\frac{1}{2}$, then for the embeddings
the corrections terms $\De_1, \De_2$ are correctly defined
as well as the phase $g_S$ is, and one has
\begin{prop} The difference of connections $A_{LC}$ and $A_0$
in the affine space $\sA_h(K^{-1}|_S)$ equals to
$$
A_{LC} - A_0 = \frac{i}{2\pi} (\De_1 + \De_2 + d {\rm ln} g_s).
\eqno (3)
$$
\end{prop}

Topologically the universal Maslov
class is the correction term for the topological type
of the covariantly constant, with respect to the determinant connection,
section of the anticanonical line bundle restricted to a lagrangian
submanifold. Suppose that a lagrangian embedding $S \subset M$ is
oriented and suppose that the restriction of the determinant
connection $A_{LC}$ to this submanifold possesses a covariant
constant section $\tilde{S} \subset P = S \times S^1$, where $P$ is the
principal $U(1)$ - bundle, associated to the restriction of the anticanonical
line bundle $K^{-1}|_S$. Then one has a simple remark

\begin{prop} The homology class of the submanifold $\tilde{S} \subset P$
has the form
$$
\begin{array}{c} [\tilde{S}] = [S] \oplus P.D.(m_S) \otimes P.D.(h) \in H_n(P, \Z)
= \\ \Z \otimes H_n(S, \Z) \oplus H_{n-1} (S, \Z) \otimes H_1(S, \Z),\\
\end{array}
$$
where the group $H_n(P, \Z)$ is decomposed according to the Kunneth
formula, and $P.D.(\al)$ denotes the homology class, Poincare
dual to the cohomology class $\al$.
\end{prop}

It follows, in particular, that in the case of  trivial
$A_{LC}$, non degenerated highest  $(n, 0)$ - form $\theta$ defines
a subset $S_1 \subset S$ via condition
$$
S_1 = \{ p \in S | \Im \theta_p = 0 \},
$$
and then the homology class $[S_1] \in H_{n-1}(S, \Z)$ is Poincare dual
to the universal Maslov class:
$$
P.D.(m_S) = [S_1].
$$

The geometrical sense of the phase $g_S$ is quite clear as well:
the dual phase transformation $\bar g_S$ sends the highest real form
$d \mu \in \Ga(\det T^*S)$ to a form of type $(n, 0)$.
Recall, that the action $\bar g_S$ on $d \mu$ is not just
the multiplication by the complex valued function $e^{- i \phi_S}$, but
indeed a gauge transformation in the corresponding $U(1)$ - bundle.
Therefore, in general we have

\begin{prop} For a lagrangian embedding $S \subset M$ to
a Kahler manifold $M$ with trivial
restriction of the determinant Levi - Civita connection the phase $g_S$
induces a real linear map
$$
g^*_S: \Ga(\det T^*S) \to \Om^{n,0}_M|_S,
$$
such that $g^*_S(f \al) = f g^*_S (\al)$ for a real function
$f \in C^{\infty}(S, \R)$.
\end{prop}

We complete the section with one more definition.

\begin{dfn} A lagrangian submanifold $S \subset M$ of a Kahler
manifold $(M, \om, I)$ is called special lagrangian if
the restriction to it of the determinant Levi - Civita connection
is trivial and the phase  $g_S$ is constant.
\end{dfn}

\section{ Examples. Stability with respect to deformations.}

{\bf Example 1.}  Let $(M = V, \om, I)$ be a symplectic vector
space with a constant symplectic form $\om$ and a constant complex
structure $I$. This means that an identification of the fibers
of the tangent bundle $TV \equiv V \times V$ is fixed and thus
a trivialization of $TV$ is fixed, defined up to
gauge $SO(n)$ - transformation and it corresponds to
a trivial connection without torsion and this connection
is exactly the Levi - Civita connection of the metric $G$,
which completes the pair $(\om, I)$ to the corresponding Kahler triple.
 The classical construction
from \cite{ArGiv}  of the universal Maslov class
of an oriented lagrangian embedding $S \subset M$ in this cases
uses the Gauss map
$$
\ga: S \to {\rm Gr}_{Lag}^{\uparrow} \times V
$$
into the grassmannization of oriented lagrangian
subspaces, and since there is the standard map
$$
\det: {\rm Gr}_{Lag}^{\uparrow} = \frac{U(n)}{SO(n)} \to U(1),
$$
then the combination $\det \cdot \ga: S \to U(1)$ gives the phase
of the lagrangian submanifold $S$ and the Maslov
class $m_S \in H^1(S, \Z)$ on it. It' s not hard to see that the classical
construction is a reduction of the construction given in the previous
section for the flat case.

Indeed, the Levi - Civita connection of metric $G$, trivializes the tangent bundle,
induces the determinant connection $A_{LC}$ on the anticanonical line bundle
$K^{-1}_M$, such that there exists a global covariantly constant with respect to
 $A_{LC}$ section $\theta^* \in
\Ga(K^{-1})$, dual to a highest holomorphic form $\theta
\in \Om^{n, 0}_V$. Since the determinant Levi - Civita
connection $A_{LC}$ is flat and  globally trivial then
after restriction to any oriented lagrangian submanifold $S \subset M$
it is automatically contained by the gauge class of the trivial
connection $A_0$. Moreover, the phase $g_S$ in this case
is related in a simple manner\footnote{
 two trivializations are related by an element from ${\rm Map} (S, U(1))$
 while two connection --- by an element from ${\rm Map}(S, U(1))/
 U(1)$, thus to pass from the first one to the second
 one needs just to add "up to constant $e^{i \cdot c}$", see
 \cite{DK}} to the gauge transformation which relates two
 trivializations: the restriction $\bar \theta |_S$ and non degenerated
 polivector field $\tau \in \Ga(\det TS)$, dual to $\pm d\mu \in
 \Ga(\det T^*S)$, where $d \mu$ is the volume form of the riemannian
 metric  $G|_S$ in an orientation (since $g_S$ is defined as a $U(1)$
 - function up to constant $e^{i \cdot c}$ then the choice
 of orientation doesn't impact on the answer ---
 due to this fact our construction doesn't depend on the orientation!).
 The last gauge transformation can be easily computed as follows:
the highest antiholomorphic form $\bar \theta$ after restriction
to $S$ can be evaluated to the highest polivector field  $\tau$
of the unit length, and since $\om$ and $I$ are compatible it follows
that for any point $p \in S$ the result
$\bar \theta_p(\tau)_p \in U(1)$. Thus we get the map
$$
\bar \theta(\tau): S \to U(1),
$$
which in principle depends on the choice of orientation
and precisely coincides with the composition $\det \cdot \ga$
from \cite{ArGiv} in the oriented case. The phase $g_S$
can be derived from $\bar \theta(\tau)$ if one forgets about
the group structure on $U(1)$, but under this the classes
$(\det \cdot \ga)^* h, g_S^* h \in H^1(S, \Z)$, where $h \in H^1(U(1) = S^1, \Z)$
is the generator, obviously coincide.

In the non orientable case for the definition of the Maslov class
one uses the following trick: instead of the map $\det$
one considers the map $\det^2$ (see \cite{ArGiv}). It's easy to see
that in terms of the construction of section 1 it corresponds
just to passage to the second power of the anticanonical line bundle
which in this case is identified to the complexification of the trivial
bundle $(\det TS)^2$; under this shift to the squares the connections
 $A_{LC}$ and $A_0$ are doubled, where the connection
$2 \cdot A_0$ is  trivial now and the gauge transformation $g_S^2$ which we
 finally get, on the one hand, is compatible with the composition
 $\det^2 \cdot \ga$, and on the other hand is the
 square of the gauge transformation $g_S$, therefore the classical
 construction from \cite{ArGiv} gives double Maslov class on the version
 of section 1.

{\bf Example 2.} It is natural to pass from the classical
flat case, presented in Example 1, to the case when for
a Kahler manifold $(M, \om, I)$ the determinant Levi - Civita connection
on the anticanonical line bundle $K^{-1}_M \to M$
admits a covariantly constant section. In this case the Ricci form
of the Kahler structure (which equals up to constant to the curvature
of the determinant connection) vanishes identically, and such a
 manifold is called the Calabi - Yau manifold\footnote{more precise,
 one has to require the simply connectedness of $M$, but according to
 a specific tradition, comes from the string theory, one calls
 "Calabi - Yau" any Kahler manifold with trivial canonical class}.
 Then for any orientable lagrangian submanifold $S \subset M$
 the restriction to it of the determinant Levi - Civita connection
 $A_{LC}$ is a flat connection with trivial periods
 and thus automatically $A_{LC} \in \sG_h(A_0)$. Therefore
in this case the correction terms  $\De_1, \De_2$ from (3)
are trivial, and for any submanifold we get
the phase $g_S$  comparing  $A_{LC}$ and $A_0$.
As it was done in Example 1, the phase $g_S$ can be
computed from the comparing of the trivializing sections for
 $A_{LC}$ and $A_0$, and the last two are dual to
 the highest holomorphic from $\theta$ (defined up to $U(1)$)
and the volume form $d\mu$ (defined up to sign) on $S$
respectively. Again one can, fixing
$\theta$, evaluate the highest holomorphic form
on the highest polivector field $\tau$ of unit length, getting the map
$$
\bar \theta(\tau): S \to U(1),
$$
which coincides up to $U(1)$ to the phase $g_S$.

The space of all lagrangian submanifolds of $M$ contains a distinguished
subset of such submanifolds $S \subset M$, that the image of the phase map $g_S$
consists of a single point:
$$
\ga_S (S) = p \in S^1.
$$
These lagrangian submanifolds are called {\it special} lagrangian;
they play the important role in the studying of the mirror symmetry effect,
see \cite{Hit}, \cite{T2}, \cite{Vafa}. From this one sees that Definition 3
of section 1 agrees with the standard terminology.

It's not hard to see that in the Calabi - Yau case
the Maslov class, defined in the previous section, is an invariant of
continuous lagrangian deformations:

\begin{prop} Let $\phi_t: S \to M$ be a lagrangian deformation
of a smooth lagrangian submanifold $S= S_0$ in a Calabi - Yau
manifold $(M, \om, I)$. Then the Maslov class $m_S$
is an invariant of this deformation:
$$
m_S = \phi_t^* m_{S_t} \in H^1(S, \Z).
$$
\end{prop}

Indeed, since for all $t$ the restriction to $S_t$ of
$A_{LC}$ belongs to the orbit $\sG_h (A_0)$, for all elements of the deformation
the phase $g_{S_t}$ is correctly defined and the topological type of
$$
g_{S_t}: S_t \to S^1
$$
doesn't depend on the deformation. But the topological type
is precisely the Maslov class $m_{S_t}$.

{\bf Example 3.} The next one in the hierarchy is
the case of Kahler - Einstein manifolds, that is
the manifolds for which
$$
\om = k \rho,
$$
where $\rho$ is the Ricci form of the Kahler metric
(this case includes the case of Calabi - Yau manifolds for which
$k = 0$). For the manifolds one has the definition of
the Maslov {\it index} for {\it Bohr
- Sommerfeld} lagrangian submanifolds, given by Fukai in \cite{Fuk}.
It's not hard to see that this Maslov index can be lifted to
a  {\it class} from $H^1(S, \Z)$, which is  the Maslov class
presented in section 1.

For a Kahler - Einstein manifold $(M, \om, I)$ the Ricci form $\rho$
is proportional to the symplectic (= Kahler) form, hence for
any lagrangian submanifold $S \subset M$ the determinant Levi - Civita
connection $A_{LC}$ is flat being restricted to $S$,
but in this case it can be non trivial. The connection $S_{LC}$ has trivial periods
 on Bohr - Sommerfeld lagrangian submanifolds only
 by the definition of these ones (see \cite{GT}). Consequently,
 for a Bohr - Sommerfeld lagrangian submanifold $S \subset M$
 the construction of section 1 does work and hence we get
 the phase $g_S: S \to S^1$. The corresponding class
 $g_S^* h, h \in H^1(S^1, \Z)$, gives the value of the Maslov index
 from \cite{Fuk}, being evaluated on the elements of $\pi_1(S)$.

However in the case of Kahler - Einstein manifolds
as well as the Maslov class is correctly defined for Bohr
- Sommerfeld lagrangian submanifolds, this class is  invariant
not for all lagrangian deformations but for a subclass consists of
 {\it isodrastic} deformations (see \cite{GT}). A lagrangian
 deformation is called isodrastic if it preserves the periods of the
 connection. In the Kahler - Einstein case for any
lagrangian embedding $S \subset M$ the correction term $\De_1$
from  (3) of the difference $A_{LC} - A_0$ is trivial, but the correction
term $\De_2$ from (3), responsible for the periods, is trivial only
for Bohr - Sommerfeld lagrangian submanifolds. If a deformation is not
isodrastic then along it the term $\De_2$ can vary on the Jacobian
$$
J_S = H^1(S, \R)/ H^1(S, \Z),
$$
in half integer point of the Jacobian the phase
is not defined at all and going around the torus one gets  the changing
by a unit of the Maslov class. But under  isodrastic deformations
the correction term $\De_2$ doesn't change and it follows that the Maslov
class is invariant. From this we have

\begin{prop} For an isodrastic deformation $\phi_t: S \to M$
of a Bohr - Sommerfeld lagrangian submanifold $S = S_0$ the Maslov class
is an invariant of the deformation
$$
m_S = \phi^*_t m_{S_t} \in H^1(S, \Z).
$$
\end{prop}

Isodrastic deformations are generated by strictly hamiltonian
vector fields and hence they are often called  hamiltonian but we use
this more sonorous term to avoid confusions which can happen by the following
reason. A vector field is called hamiltonian if it preserves the symplectic
form. But this vector field induces the deformation which can be non isodrastic
if the field is not generated by a function that is not strictly hamiltonian.

At the end, in the general situation
a lagrangian deformation can be collected by lagrangian embeddings
with bounded periods (for each element of such a deformation
the monodromy of the connection $A^1_{LC}$ is  equal to
-1 for no primitive elements from $\pi_1(S)$), and for the elements
of this deformation

--- the correction terms $\De_2$ are correctly defined;

--- the phases $\phi_{S_t}$ are correctly defined;

--- the Maslov class is correctly defined and is an invariant of
the deformation family.

It can be formulate, f.e., as follows:

\begin{prop} Let $\phi_t: S \to M$ be a lagrangian deformation
of $S = S_0$ with bounded periods. Then the Maslov class $m_{S_t}$
is correctly defined via the correction connection $A^0_{LC}$
and is an invariant of the deformation:
$$
m_S = \phi_t^*(m_{S_t}) \in H^1(S, \Z).
$$
\end{prop}

The last proposition is just an illustration
how the construction of section 1 can be exploited in more general situation
than Bohr - Sommerfeld lagrangian embeddings to Kahler - Einstein
manifolds. However the last case is of our main interest
in view of possible applications to algebraic lagrangian geometry.

\section{Minimality}

The  problem of  minimality for riemannian volume
of lagrangian submanifolds and the possibility of the deformation
to a minimal lagrangian submanifold are
solved long ago for the case of symplectic vector spaces, see \cite{Mor}.
The compact case has been studied in \cite{Br} and \cite{Oh1},
\cite{Oh2} --- first, Bryant proved that a lagrangian submanifold
$S \subset M$ is minimal only if the restriction of the Ricci form
$\rho$ to $S$ is trivial and the determinant Levi - Civita connection
$A_{LC}$ admits a covariantly constant section that is its periods are trivial;
then Oh reduced the minimality problem to the consideration of the
Hodge decomposition of the  {\it mean curvature} of the lagrangian
submanifold.

Recall some definitions from riemannian geometry. Let $S
\subset M$ be an embedding to a riemannian manifold. Then it is defined the
 {\it second quadratic form}
$$
{\rm II}: TS \to N \otimes T^*S,
$$
where $N$ is the normal bundle. Namely, sections of $TS$ are differentiated
as sections of $TM|_S$ with respect to the total Levi - Civita
connection $\nabla_{LC}$ on the  "big" tangent bundle $TM|_S$,
and then the result is projected to the normal component in
$TM|_S \otimes T^*S$. It's not hard to see that this composition
is a tensor --- not a differential operator ---
and this tensor is called the second quadratic form. At the same time
the total Levi - Civita connection $\nabla_{LC}$ on $TM|_S$
can be recovered from the reduced Levi - Civita connection $\nabla_{LC}^S$,
induced by the restriction of the riemannian metric $G$ to  $S$,
and the second quadratic form ${\rm II}$.

Further, in our situation $S \subset M$ is a lagrangian embedding
to a symplectic manifold $(M, \om)$, and the choice
of a compatible riemannian metric $G$ attaches to $S$
certain real 1- form $\al_H$,  called the
{\it mean curvature}, by the following rules: the trace of the second
quadratic form
$${\rm tr \quad II}
$$
is a section of the normal bundle
and since for a lagrangian embedding the normal bundle  $N$
is isomorphic to $T^*S$, see \cite{Wein},
 then the trace is represented by a section of $T^*S$, that is by a  1 -form.

 According to a Hodge theory, any form $\al$ in presence of a metric
 is decomposed into three parts
 $$
 \al = \al_1 + \al_0 + \al_{-1},
 $$
 where $\al_{ \pm 1}$ is  (co) exact and $\al_0$ is harmonic,
 and this decomposition is unique (see \cite{DK}). We have
 $$
 \al_1 + \al_0 \in \ker d, \quad \al_0 + \al_{-1} \in \ker d^*
 $$
 and
 $$
 \al_0 \in \ker d \cap \ker d^* = \sH^1,
 $$
 where the last intersection is exactly the space of harmonic forms.

 It's known (see \cite{Br}, \cite{Oh1}), that a lagrangian submanifold
 $S_0$ is minimal with respect to local lagrangian deformations
 (L - minimal) if and only if its mean curvature $\al_H$ is trivial.
 In \cite{Oh1} one proposes the notion of  H- minimality ---
 the minimality with respect to isodrastic deformations ---
 and one shows that a lagrangian embedding is H - minimal
 if and only if $d^* \al_H = 0$. Moreover, Oh proved that
 for a lagrangian embedding to a Calabi - Yau manifold the mean curvature
 is always closed and  represents an integer cohomology class.

 The constructions of section 1 are related to the minimality
 problem by the following identity.

 \begin{prop} Let $S \subset M$ be an orientable lagrangian
 embedding to a Kahler manifold $(M, \om, I)$. Then the connections
  $A_{LC}, A_0$, defined in section 1, and  1 - form $\al_H$
 are related by identity:
 $$
2 \pi i \al_H =  (A_{LC} - A_0).
\eqno (4)
 $$
 \end{prop}

 As a corollary we get, in particular,
 the coincidence of the restricted Ricci form and the differential
 of the mean curvature,
  $$
 \rho|_S = d \al_H,
 $$
which has already been established in \cite{Daz},
since
 $$
 \rho|_S =  - \frac{i}{2\pi} F_{LC} = \frac{1}{2\pi i} d(A_{LC} - A_0).
 $$

 Let us stress that such an identity is possible in the case
 of integrable complex structures only.

The proof of Proposition 9 is based on the reconstruction
of the total Levi - Civita connection $\nabla_{LC}$ on $TM|_S$
from the reduced connection $\nabla_{LC}^S$ on $TS$
and the second quadratic form ${\rm II}$.

 Indeed, let's represent $TM|_S$ as the complexification of $TS$
 and consider the holomorphic tangent bundle $T^{1,0}M|_S$,
 which is a $U(n)$- bundle. Then it carries two
 $U(n)$ - connections: the total Levi - Civita connection $\nabla_{LC}$
 (by the intergability condition the holonomy of
 metric $G$ is $U(n)$) and the connection, which comes
 after the complexification from the reduced Levi - Civita connection
 and we will denote this one by the same symbol
$\nabla_{LC}^S$ to simplify the explanation. Any two $U(n)$- connections
over a  $U(n)$ - bundle are differ by a  1- form with values in the
Lie algebra $u(n)$. It's not hard to see that for our $U(n)$ - connections
$\nabla_{LC}$ and $\nabla_{LC}^S$ this  1- form is given
by the second quadratic form being considered as a 1- form with values
in symmetric endomorphisms of $TS$. It is very simple to
convert a symmetric endomorphism of  $TS$ to the corresponding
skew symmetric endomorphism of the complexification $TS \otimes U(1)$,
and thus one sees that
  $$
  \nabla_{LC} - \nabla_{LC}^S = 2 \pi i {\rm II}.
\eqno (4')
  $$
Further, by the definition of our connections $A_{LC}$ and $A_0$
on the restriction on the anticanonical line bundle
$K^{-1}|_S = \det T^{1,0} M|S$ one gets, reducing the situation in (4')
to the determinants, that the identity holds
$$
A_{LC} - A_0 = 2 \pi i {\rm tr} {\rm II},
$$
and thus taking into account the definition of the mean curvature
it follows the statement of Proposition 9.

Thus for any lagrangian embedding $S \subset M$ to a Kahler
manifold $(M, \om, I)$ there is the same 1-form, the mean curvature $\al_H$,
which can be in view of identity (4) decomposed into components
in two ways: as a 1 - form in the Hodge theory and as
a connection form on a trivial bundle in the theory of
connections:
$$
\begin{array}{l} \De_1 + \De_2 + d \ln g_S = \al_H = \\
 (\al_H)_1 + (\al_H)_h + (\al_H)_{-1}.\\
 \end{array}
$$
The components of the decompositions are related as follows:
$$
\begin{array}{l}
\De_2 + d \ln g_S = (\al_H)_1 + (\al_H)_h \in \ker d\\
\De_1 = (\al_H)_{-1} \in \Im d^*,\\
\end{array}
$$
and on the real cohomology level
$$
[(\al_H)_h] = [\De_2] + m_S \in H^1(S, \R).
$$
Therefore in the generic case the Maslov class is the "integer part"
of the cohomology class which is presented by the harmonic part
of the mean curvature.

The identity (4) implies  a number of corollaries.

\begin{cor} A lagrangian submanifold $S \subset M$ is L - minimal
if and only if the restriction to it of the anitcanonical line bundle
with the determinant Levi - Civita connection admits a covariantly constant
section and the phase $g_S$ is constant.
\end{cor}

Indeed, the existence of a covariantly constant section is equivalent to
$\De_1 = \De_2 = 0$ in the decomposition (3), and the rest component
$d \ln g_S$ is trivial if and only if the phase is constant. In other words
a minimal lagrangian submanifold must be  special lagrangian.

\begin{cor} A lagrangian submanifold $S \subset M$ with trivial
restriction of the determinant Levi - Civita connection and
trivial Maslov class is H- minimal if and only if it is
L - minimal or special lagrangian.
\end{cor}

Really according to  \cite{Oh1}, the H - minimality is equivalent
to the condition $d^* \al_H = 0$;
the components $\De_1, \De_2$ from (3) always lie in the kernel of $d^*$,
and if the Maslov class $m_S$ is trivial then the component
$d \ln g_S$ is an exact form and it lies in the kernel of $d^*$
if and only if it vanishes.

Finally, the propositions of section 2 together with the identity (4) show
that

\begin{cor} 1. A lagrangian submanifold $S \subset M$ of a Calabi - Yau
manifold $M$ with nontrivial Maslov class can not be
transported to a minimal lagrangian submanifold by lagrangian
deformations.

2. A Bohr - Sommerfeld lagrangian submanifold $S \subset M$ of
a Kahler - Einstein manifold $M$ with nontrivial
Maslov class can not be transported to a minimal
lagrangian submanifold by isodrastic lagrangian deformations.
\end{cor}

General case requires more detailed considerations
and since the applications we need are contained by the case
of Kahler - Einstein manifolds we stop here with the presented
relationships.

\section{Applications}

The definition of the phase $g_S$ of a lagrangian embedding
to a Kahler manifold given at section 1 of the present paper
makes it possible to develop some approaches to the problems
of mirror symmetry and geometric quantization
(see \cite{T1}, \cite{Vafa}, \cite{GT}).

{\bf Special lagrangian geometry}. There is an approach to the
mirror symmetry problem which is based on the studies of the moduli
spaces of special lagrangian submanifolds of a Calabi - Yau
manifold (see \cite{T2}, \cite{Vafa}).
Special lagrangian submanifolds are mirror
dual to stable holomorphic bundles on the mirror partner, but
the theory of the moduli spaces
of special lagrangian submanifolds of Calabi - Yau
manifolds, presented in \cite{Joyce}, is not completely
finished yet, and the main difficulties come from
the singularities which limiting special lagrangian submanifold
can carry. A.N. Tyurin some time ago proposed to study even more
restricted class --- the class of Bohr - Sommerfeld special lagrangian
submanifolds whose moduli space has virtual dimension zero.
The same question has sense for the case of Kahler - Einstein manifolds as well
and we consider this case in this section.

If $(M, \om, I)$ is a Kahler - Einstein manifold then for each orientable lagrangian
submanifold $S$ the restriction of the determinant Levi - Civita
connection $A_{LC}$ is a flat connection. By the definition  of Bohr
- Sommerfeld lagrangian submanifolds, see \cite{GT},
$$
\sM_{SpBS} = \sM_{BS} \cap \sM_{SpLag}
$$
consists of lagrangian submanifolds with trivial
period part $\De_2$, for which the phase $g_S$ and the Maslov class $m_S$
are correctly defined and moreover the first one is constant
and the last one is trivial. Codimension of $\sM_{BS}$
 in the space of all lagrangian submanifolds
 equals to $b_1(S)$ (see \cite{GT});
 due to Corollary 2 the dimension $\sM_{SpLag}$ coincides with the dimension
 of the space of H- minimal lagrangian submanifolds which is
 computed by Oh in general case in the paper \cite{Oh2}
 and is equal to $b_1(S)$. The H - minimality is quite reasonable
 to investigate in connection with the Bohr - Sommerfeld condition
 since locally the moduli space  $\sM_{BS}$ is exactly
 generated by isodrastic deformations, see \cite{jT1}. Therefore
 one could expect that the dimension of $\sM_{SpBS}$ is zero:
 $$
 \sM_{SpBS} = \{ p_1, ..., p_d\}
 $$
 and the number of the points, the space $\sM_{SpBS}$ consists of,
 is a symplectic invariant. However in the same paper
 Oh presented an example when the real dimension of the space of
 H- minimal lagrangian submanifold is greater than the virtual one.
 Briefly recall the example and discuss why this jumping
 takes  place.

Let  $(M, \om, I)$ be  projective line $\C \PP^1$ endowed with the
standard Fubini - Study metric. Then (see the toy example
from the last section of \cite{GT}) a smooth loop $\ga \subset M$ (
which is always lagrangian due to the dimensional reason) is
Bohr - Sommerfeld if and only if it divides the surface of $S^1 = \C \PP^1$
into two parts with the same area. On the other hand,
a smooth loop $\ga$ is minimal  (H - minimal) if and only if
it is a big circle (a circle). Therefore the local dimension
of $\sM_{\min}$ equals to 2, while the virtual dimension
equals to $b_1(S) = 1$. Respectively, the intersection
$$
\sM_{BS} \cap \sM_{\min} = \sM_{SpBS}
$$
has dimension 2, while we expect it equals to zero. In \cite{Oh2} one explains
this fact by the existence of a special symmetry of the Fubini - Study metric;
if we deform this structure to
such a structure that the sphere is transformed to  a ball for the american football
then for this Kahler structure the space $\sM_{\min}$ is already
1 dimensional and the intersection
$$
\sM_{BS} \cap \sM_{\min} = \sM_{SpLag}
$$
in this case consists of the following two components --- a single point
(which corresponds to equator) and a 1 -dimensional component (which
corresponds to meridians). Deforming the Kahler structure to even more
generic one we would get the required zero dimensional space $\sM_{SpLag}$.

The cause of the dimensional jumping of the space $\sM_{SpLag}$ is
hidden in the existence of some infinitesimal symmetries of the Kahler structure
which correspond to {\it quasisymbols} on $(M, \om, I)$, that is
the smooth real functions whose hamiltonian vector fields preserve
whole the Kahler structure (all details see in \cite{jT1}).
For the Fubini - Study Kahler structure on $S^2$ the quasisymbol
space has dimension 2; for the "american ball" the dimension
of the quasisymbol space equals to 1 and in this case the equator
is invariant with respect to their action while the meridians are not
and it follows that we have two components of different dimensions.

Another example, presented in \cite{jT2}: one takes as
$(M, \om, I)$ an elliptic curve and since
an elliptic curve doesn't admit quasisymbols the space
$\sM_{SpLag}$ is zero dimensional (and finite).
Quasisymbols exist neither for a "real" Calabi - Yau
manifold (it follows from the simply connectedness)
no for an abelian variety. It was conjected that for a
Kahler - Einstein manifold
 $(M, \om, I)$, which doesn't admit quasisymbols,
 the space $\sM_{SpLaG}$ is zero dimensional and
 the number of points, which form $\sM_{SpLag}$,
 is a symplectic invariant. In general case the dimension of
 $\sM_{SpLag}$ coincides with the dimension of the quasisymbol
 space, which acts on $\sM_{SpLag}$, and the number of connected
 components of $\sM_{SpLag}$ is a symplectic invariant.

 To prove the conjecture it suffices
 to express the phase for the deformation of a Bohr - Sommerfeld
 special submanifold $S_0 \in \sM_{SpLag}$,
 induced by a smooth function $f \in C^{\infty}(S_0, \R)$
 (all details of the correspondence between functions
 and isodrastic deformations see in \cite{GT}).
 Since  $S_0$ has trivial Maslov class, being
 special lagrangian, any isodrastic deformation of $S_0$
 has trivial Maslov class too and hence the phase
 of $S_t$ is presented by the form
 $$
 e^{i \phi_S},
 \eqno (5)
 $$
 where $\phi_S$ is a smooth function,
defined up to an additive constant (precisely as  $f$ is, see \cite{GT}).
 A naive suggestion that $f = \phi$, of course,
 is false ---  the functions are related by a differential equation
 with the principal term, equals to the Laplace operator.
 The precise formula, which relates the phase with the deformation,
 would lead not only to checking of the conjecture, presented above,
 but as well to the definition of some half weighting rules
 for Bohr - Sommerfeld lagrangian submanifolds.

 {\bf Weighting and half weighting rules.}

The moduli space of half weighted Bohr - Sommerfeld
lagrangian cycles introduced in \cite{GT},
plays an important role in quantization
of classical mechanical systems, see \cite{jT1}.
Consider a Kahler - Einstein manifold
$(M, \om, I)$ and suppose that $[\om] \in H^2(M, \Z)$ ---
that is the prequantization condition holds. Then (see \cite{GT})
one chooses prequantization data $(L, a)$, and in our case when
$$
[\om] = k \cdot K_M
$$
the prequantization bundle can be chosen
together with an isomorphism, identified the hermitian structures
on $L$ and $K_M$. Then the prequantization connection $a \in sA_h(L)$
can be chosen such that $a_K \in \sA_h(L^k) = \sA_h(K_M)$
is contained by the same gauge class as the connection ${\tilde A}_{LC}^*$
is, where ${\tilde A}_{LC}$ is the determinant Levi - Civita connection
on whole $M$. Let's fix a positive volume $r > 0$
and consider the moduli space of half
weighted Bohr - Sommerfeld lagrangian cycles $\sB_{BS}^{hw,r}$
of the fixed volume (the definition see in \cite{GT}).
Let's take the component $\sB_{BS,0}^{hw,r}
\subset \sB_{BS}^{hw,r}$, which consists of lagrangian
submanifolds with trivial Maslov class. Then we have a bi - section
$$
\sB_{BS, 0} \to \sB^{hw,r}_{BS,0},
$$
defined by the phase half weighting rule
$$
\theta_S^2 = (\phi_S + c) d \mu,
$$
where $\phi_S$ is defined in (5), $ d \mu$ is the volume form
induced by the riemannian metric $G$ on $S$,
and constant $c \in \R$ is defined by the condition
$$
\int_S (\phi_S + c) d \mu = r.
$$
It's easy to see that under this rule minimal lagrangian submanifolds
are distinguished by the fact that for each minimal one
the square of the canonical half weight $\theta_S^2$ is proportional
to the riemannian volume form.
$$
$$
At the end we remark that it is possible
to give an interpretation for the Maslov class,
presented in section 1 of this paper,
in non integrable case as well. But we postpone
the discussion of this question until the time when some reasonable
geometric or topological applications will be found hoping
that this time will come soon.

\end{document}